\documentclass[1p,final]{elsarticle}

\makeatletter
\def\ps@pprintTitle{%
 \let\@oddhead\@empty
 \let\@evenhead\@empty
 \def\@oddfoot{}%
 \let\@evenfoot\@oddfoot}
\makeatother

\usepackage{amssymb}
\usepackage{amsmath}
\usepackage{url}

\begin{document}

\begin{frontmatter}

\title{Using incomplete indefinite $LDL^T$ preconditioning for inexact
interior point methods for linear programming\tnoteref{t1}}
\tnotetext[t1]{APMOD 2012 extended abstract}

\author{Robert Luce}

\date{November 2011}

\ead{luce@math.tu-berlin.de}

\ead[url]{http://www.math.tu-berlin.de/\~{}luce}

\address{
Technische Universit\"at Berlin\\
Institut f\"ur Mathematik, MA 3-3\\
Stra{\ss}e des 17. Juni 136\\
10623 Berlin, Germany}

\begin{abstract}
Most linear algebra kernels in interior point methods for linear
programming require the solution of linear systems of equation with
the matrix $N = A^TD^{-1}A$ (or $AD^{-1}A^T$), where $A$ denotes the
constraint matrix of the linear program.  This matrix $N$ arises from
the reduced KKT system by block elimination.  If the number of
non-zeros in $N$ or in its Cholesky factorization $N= LL^T$ is very
large, the computational cost and memory requirement to solve the
linear systems of equations with $N$ may be prohibitively large.  In
this work we implement an interior point method described by R. Freund
and F. Jarre~\cite{FreundJarre:1996}.  Forming the normal equation
matrix $N$ is avoided altogether and we work with the reduced KKT
system instead.  We solve the linear systems for the Newton directions
iteratively only to low accuracy using SQMR and an indefinite
multilevel preconditioner
(PARDISO~\cite{SchenkBolle:2008,HagemannSchenk:2007}).  Preliminary
numerical results are encouraging.
\end{abstract}

\begin{keyword}
    \MSC 90C05 \sep 65K05 \sep 65F10 \sep 65F50
\end{keyword}

\end{frontmatter}


\section{Introduction}
\label{sec:intro}
We investigate an iterative alternative for the linear algebra kernel
in interior point methods (IPM) for linear optimization. By that we
understand that the linear systems of equations that arise in the
course of the IPM are solved only to a certain (low) relative accuracy
using an iterative solution method.  A typical situation where the
usual direct approach for solving those systems using a Cholesky
factorization of the normal equations is not satisfying arises when
the normal equations or, more often, its Cholesky factorization
suffers a large amount of fill.  So the setting we have in mind is the
case where forming the normal equation matrix or its Cholesky
factorization is not desired and thus we work exclusively with the
augmented indefinite KKT system.\footnote{Another situation where the
augmented system {\itshape may} be preferred is the presence of free
variables in the primal LP, which are conceptually simpler to treat in
the augmented system.} We note in passing the computational cost of an
matrix vector product (the only relevant cost in iterative linear
systems solvers) is the same for the augmented system and the normal
equation matrix.  Further it can easily be seen that every
preconditioner for the normal equation matrix is mathematically
equivalent to one for the augmented KKT matrix (see, for example,
\cite{OliveiraSorensen:2005}) but not vice versa.

A large number of publications on this subject and, more generally,
inexact Newton methods for nonlinear systems of equations exist.  We
do not attempt to cite all the relevant ones in this context.  We
remark that we do not assume any knowledge of structure in the LP
constraint matrix (in case such structure exists) or take advantage of
heuristics that try to identify an optimal basis partitioning as in
\cite{OliveiraSorensen:2005, GondzioHall:2008}.   

Our work was prompted by promising numerical results obtained in
non-convex optimization \cite{SchenkWaechter:2008} and a certain
eigenvalue problem \cite{SchenkBolle:2008} where an indefinite
multilevel preconditioner is used that is part of the software package
PARDISO\footnote{\url{http://www.pardiso-project.org}}.
Further our work is based on the inexact IPM described by R. Freund
and F. Jarre~\cite{FreundJarre:1996}.

\section{Statement of the problem}

We consider the primal-dual linear programming pair
\begin{align*}
\mbox{min}\;  & c^T x   & \mbox{max}\;  & b^T y \\
\mbox{s.t.}\; & Ax = b  & \mbox{s.t.}\; & A^T y + z = c\\
            & x \ge 0   &               & z \ge 0.
\end{align*}
For a given iterate $(x,y,z)$ and a centering parameter $\mu$ the
residuals for the corresponding first order optimality system are
\begin{align*}
r_p = Ax - b, \quad r_d = A^T y + z - c, \quad r_c = Xz - \mu e.
\end{align*}
(For $x \in \mathbb{R}^n$ the diagonal matrix $\text{diag}(x)$ is
denoted by $X$.) Then the Newton direction $[d_y,d_x,d_z]^T$ at this
particular step is implied by the solution of
\begin{equation}
\label{eqn:augmented_KKT}
    K \tilde{d} =
    \begin{bmatrix}
        0   & A\\
        A^T & X^{-1}Z
    \end{bmatrix}
    \begin{bmatrix}
        d_y\\
        \tilde{d}_x
    \end{bmatrix}
    = -
    \begin{bmatrix}
        r_p\\
        r_d - X^{-1}Z r_c
    \end{bmatrix}
    = \tilde{r},
\end{equation}
where $d_x, d_z$ can be recovered from $\tilde{d}_x$ by block
elimination with $X$ and $Z$.  Our goal is to solve
\eqref{eqn:augmented_KKT} up to a certain accuracy $\eta \in (0,1)$,
so that $\|\tilde{r} - K \tilde{d}\| \le \eta \|\tilde{r}\|$.  Here
$\eta$ has to be chosen such that the norms of the residuals $r^+_p$
and $r^+_d$ of the next iterate are sufficiently smaller than the
current ones (As pointed out in \cite{FreundJarre:1996}, $\|r^+_c\|
\approx 0$ can always be achieved.).  The solution of
\eqref{eqn:augmented_KKT} should require significantly less memory
than the solution of the corresponding normal equation.

\section{The PARDISO multilevel preconditioner}

We now describe very briefly the techniques in the PARDISO multilevel
preconditioner; details can be found in \cite{SchenkBolle:2008}.  The
preconditioner is based on an incomplete multilevel indefinite $LDL^T$
factorization, where a specified bound $\kappa$ on the norm of
$\|L^{-1}\|$ is guaranteed.  In the course of the factorization, all
pivots (1-by-1 or 2-by-2) that fail to satisfy the guaranteed bound on
the norm of $L^{-1}$ are postponed and constitute the next level of
the factorization after the Schur complement has been built.  The
dropping of entries happens both when forming the Schur complements or
as a means to keep the resulting $L$ factor sparse.  Both dropping
tolerances can be controlled independently.  The factorization is
preceded by a reordering strategy that aims at reducing the resulting
fill in the $L$ factor as well as increasing the (block) diagonal
dominance of the reordered matrix.

The iterative solver that is built-in in PARDISO is based on the
simplified QMR iteration (see \cite{FreundJarre:1996} and references
within).  Part of this work is the discussion and comparison of
various heuristics that control the parameters of the preconditioner
in the course of the IPM iteration as well as the termination
criteria for the QMR iteration.

\section{A motivating example}

As outlined in the introduction we focus on LP instances where the
fill in the Cholesky factor of the normal equation matrix is large.  A
set of LP instances with that characteristic was provided to us by the
Zuse Institute Berlin, from which the example in table
\ref{tab:example} is taken.

\begin{table}
\begin{center}
\begin{tabular}{lll}
{} & CPLEX & inexact IPM\\ \hline
IPM iterations & 17   & 67\\
fill ratio (w.r.t. $A$) & 210 & 4.5 (avg.)\\
fill ratio (w.r.t. $AA^T$) & 42 & n/a\\
solution time         & $\sim$18,000s & $\sim$3,200s
\end{tabular}
\caption{Comparison for an LP instance with 248.127 columns,
386,294 rows and 3,254,730 non-zeros in $A$.}
\label{tab:example}
\end{center}
\end{table}

We compare our implementation with the interior point method that
comes with CPLEX (we used version 12.1).  The linear algebra kernel in
that method is based on the Cholesky factorization of the normal
equations.  In the example above the fill ratios shown refer to the
quotient of the number of nonzeros in the Cholesky factor of the
normal equation matrix and the number of nonzeros in $A$ (lower
triangle of $AA^T$, respectively).  Note that storing the Cholesky
factor requires about 6GB of memory, while the preconditioner requires
much less memory.  The total solution time for CPLEX is about a factor
of five larger than the total solution time of the inexact IPM.

In our work we present a detailed numerical comparison of the
performance of our inexact IPM implementation with CPLEX on linear
programming instances that suffer much fill in the Cholesky factor.

\bibliographystyle{elsarticle-num}
\bibliography{extabs_luce}

\begin{thebibliography}{1}
\expandafter\ifx\csname url\endcsname\relax
  \def\url#1{\texttt{#1}}\fi
\expandafter\ifx\csname urlprefix\endcsname\relax\def\urlprefix{URL }\fi
\expandafter\ifx\csname href\endcsname\relax
  \def\href#1#2{#2} \def\path#1{#1}\fi

\bibitem{FreundJarre:1996}
R.~W. Freund, F.~Jarre, A {QMR}-based interior-point algorithm for solving
  linear programs, Math. Programming 76~(1, Ser. B) (1996) 183--210, interior
  point methods in theory and practice (Iowa City, IA, 1994).

\bibitem{SchenkBolle:2008}
O.~Schenk, M.~Bollh{\"o}fer, R.~A. R{\"o}mer, On large-scale diagonalization
  techniques for the {A}nderson model of localization, SIAM Rev. 50~(1) (2008)
  91--112.

\bibitem{HagemannSchenk:2007}
O.~Schenk, A.~W{\"a}chter, M.~Hagemann, Matching-based preprocessing algorithms
  to the solution of saddle-point problems in large-scale nonconvex
  interior-point optimization, Comput. Optim. Appl. 36~(2-3) (2007) 321--341.

\bibitem{OliveiraSorensen:2005}
A.~R.~L. Oliveira, D.~C. Sorensen, A new class of preconditioners for
  large-scale linear systems from interior point methods for linear
  programming, Linear Algebra and its Applications 394 (2005) 1--24.

\bibitem{GondzioHall:2008}
G.~Al-Jeiroudi, J.~Gondzio, J.~Hall, Preconditioning indefinite systems in
  interior point methods for large scale linear optimisation, Optim. Methods
  Softw. 23~(3) (2008) 345--363.

\bibitem{SchenkWaechter:2008}
O.~Schenk, A.~W{\"a}chter, M.~Weiser, Inertia-revealing preconditioning for
  large-scale nonconvex constrained optimization, SIAM J. Sci. Comput. 31~(2)
  (2008/09) 939--960.

\end{thebibliography}

\end{document}